\documentclass[11pt,showlabels]{article1}
\usepackage{amssymb}
\usepackage{epsfig}
\usepackage{amsmath}
\usepackage{graphicx}
\usepackage{subfig}
\usepackage{appendix}
\usepackage{colortbl}
\usepackage{color}
\newtheorem{problem}{Problem}
\newtheorem{theorem}{Theorem}%[section]
\newtheorem{assumption}{Assumption}%[section]
\newtheorem{lemma}{Lemma}%[section]
\newtheorem{remark}{Remark}%[section]
\newtheorem{corollary}{Corollary}%[section]

\title{\bf  Deterministic Optimal Control of It$\hat{o}$ Stochastic Systems with Random Coefficients
 \thanks{
 This work is supported by the National Natural Science
Foundation of China under Grants  61573221,
61633014.
$^{*}$Corresponding author: Huanshui Zhang. Email: hszhang@sdu.edu.cn}
}

\author{Hongdan Li,\ Juanjuan Xu and\ Huanshui Zhang$^{*}$
\ \\
\\
\ \ \  School of Control Science and Engineering, Shandong
University, \\Jinan Shandong 250061, China}

\begin{document}

\pagenumbering{arabic}
 \setcounter{page}{1}

\pagenumbering{arabic} \thispagestyle{empty} \setcounter{page}{1}

%{\footnotesize \noindent

%\vspace{-20pt}
%\begin{center}\section{Introduction}
\baselineskip 16pt
\date{}
%\begin{titlepage}
  \maketitle
\begin{abstract}
This paper is concerned with the deterministic optimal control of It$\hat{o}$ stochastic systems with random coefficients. The necessary and sufficient conditions for the unique solvability of the optimal control problem with random coefficients are derived via the solution to the coupled stochastic Riccati-type equations. An explicit expression of the deterministic optimal controller for this problem is given. The presented results include the case  of  deterministic coefficient \cite{15} as special case.
\bigskip

\noindent \textbf{Keywords:} Deterministic optimal control; random coefficients;  FBSDEs;  the coupled stochastic Riccati-type equations.
%\end{titlepage}
\end{abstract}

\pagestyle{plain} \setcounter{page}{1}
\section{Introduction}
The linear quadratic (LQ, for short) control problem which was first discussed by Kalman \cite{1} is one of the classical yet important problem in both theory and engineering applications and has received great attentions since 60's of last century. Initially, the deterministic systems were considered; see \cite{2}-\cite{4} and references therein.

The LQ optimal control problem was extended to stochastic system with multiplicative noise due to the practical applications; see \cite{6}-\cite{15} and references therein.  Stochastic LQ problem was pioneering studied by Wonham \cite{5}.  By using a generalized Riccati equation, J. Bismut solved the standard stochastic LQ  problem in \cite{11}, and the result was latter extended to indefinite stochastic LQ problem in \cite{6}. In \cite{13} and \cite{14}, Tang considered the general linear quadratic optimal stochastic control problem with random coefficients via the stochastic maximum principle and the dynamic programming principle, respectively. And the mixed optimal control of a linear stochastic system with a quadratic cost functional was studied in \cite{15}.

Different from the previous works, in this paper, we aim to find a deterministic controller for the It$\hat{o}$ stochastic systems with random coefficients. As compared with previous work \cite{15}, the problem is much involved due to the random coefficient. Actually,  the optimal controller of this problem is not simply a conventional feedback form $K_{t}x(t)$ as in \cite{13} and \cite{14} or mean-field type $K_{t}E[x(t)]$ as in \cite{15}. It remains challenging to derive an explicit expression of the deterministic optimal controller  due to the correlation between the state and coefficients.

Motivated by \cite{16} and \cite{17} in which the substantial progress for the optimal LQ control has been made by solving  forward and backward stochasitic differential/difference equations (FBSDEs), this paper focuses on the deterministic optimal control  of It$\hat{o}$ stochastic systems with random coefficients.  An analytical solution to the FBSDEs associated with the control problem is derived and then the necessary and sufficient condition for the unique solvability of linear quadratic optimal control problem with random coefficients is given via a coupled Riccati-type equation. The deterministic optimal controller is accordingly designed with Riccati-type equation.

The rest of this article is structured as follows. Section 2 gives the problem statement. The main result of this paper is presented in Section 3.   A summary is shown in Section 4.

\emph{Notation }: \ $\mathcal{R}^{n}$ is the $n$-dimensional Euclidean space and $\mathcal{R}^{m\times n}$  the norm bounded linear space of all
$m\times n$ matrices. $Y'$ is  the transposition of $Y$ and  $Y\geq 0 (Y>0)$ means  that  $Y\in \mathbb{R}^{n\times n}$ is symmetric positive semi-definite  (positive definite). Let $(\Omega,\mathcal{F},P,\{\mathcal{F}_{t}\}_{ t\geq0} )$ be a complete  probability space with natural filtration
$\{\mathcal{F}_{t}\}_{ t\geq0}$ generated by the standard Brownian motion $w(t)$ and system initial state augmented by all the $\mathcal{P}$-null sets.

%%%%%%%%%%%%%%%%%%%%%%%%%%%%%%%%%%%%%%%%%%%%%%%%%%%%%%%%%%%%%%%%%%%%%%%%%%%%%%%%%%%%%%%%%%%%%%%%%%%%%%%%%%%%%%%%%%%%%%%%%%%%%%%%%%%%%%%%%%%%%%%%%%%%%%%%%%%%%%
\section{Problem  Statement  }

\normalsize

Consider the It$\hat{o}$ stochastic system
\begin{eqnarray}
\left\{
\begin{array}{lll}
dx(t)=[A_{t}x(t)+B_{t}u(t)]dt+[C_{t}x(t)+D_{t}u(t)]dw(t), \\
 x(0)=x_{0},
\end{array}
\right.   \label{f1}
\end{eqnarray}
where $x(t) \in \mathcal{R}^{n} $ is the state; $u(t) \in \mathcal{R}^{m}$ is the control input; $w(t)$ is the one-dimensional standard Brownian motion. $x_{0} \in \mathcal{R}^{n}$ is the initial value.  \\
 The cost functional is given as following
\begin{eqnarray}
J_{T}=E\Big\{\int_{0}^{T}\left[x'(t)Q_{t} x(t)+u'(t)R_{t} u(t)\right]dt+x'(T)P_{T} x(T)\Big\}.\label{f2}
\end{eqnarray}
First, we will make the following assumption.
\begin{assumption}
The matrix processes $A_{t},C_{t}, Q_{t}:[0,T]\times\Omega\rightarrow \mathcal{R}^{n\times n}$, $B_{t},D_{t}:[0,T]\times\Omega\rightarrow \mathcal{R}^{n\times m}$, $R_{t}:[0,T]\times\Omega\rightarrow \mathcal{R}^{m\times m}$ and the random matrix $P_{T}:\Omega\rightarrow \mathcal{R}^{n\times n}$ are uniformly bounded and $\mathcal{F}_{t}$-adapted or $\mathcal{F}_{T}$-measurable. $Q_{t}, R_{t}, P_{T}$  are positive semi-definite matrices.
\end{assumption}
The problem to be considered is formulated as
\begin{problem}
Find a deterministic optimal control $u(t)$  to minimize (\ref{f2}) subject to  (\ref{f1}).
\end{problem}
\begin{remark}
Different from the previous works, such as \cite{13} and \cite{14} in which the purpose is to find a $\mathcal{F}_{t}$-adapted optimal controller $u(t)$ to minimize the cost functional, or \cite{15} in which the coefficients throughout the paper are all required to be deterministic, this paper will discuss the deterministic optimal control of It$\hat{o}$ stochastic systems with random coefficients.
\end{remark}

\section{Main Result}

In order to find the solution to Problem 1, we will introduce the following lemma first, i.e., the stochastic maximum principle.
\begin{lemma}
 Let Assumption 1 be satisfied. Problem 1 is uniquely solvable if and only if the following FBSDEs have a unique solution.
\begin{eqnarray}
\left\{
\begin{array}{lll}
0=E[B_{t}'p(t)+R_{t}u(t)+D_{t}'q(t)],\\
dp(t)=-[A_{t}'p(t)+C_{t}'q(t)+Q_{t}x(t)]dt+q(t)dw(t), \\
 dx(t)=[A_{t}x(t)+B_{t}u(t)]dt+[C_{t}x(t)+D_{t}u(t)]dw(t), \\
 x(0)=x_{0},p(T)=P(T)x(T).\label{f3}
 \end{array}
\right.
\end{eqnarray}
 \end{lemma}
\emph{Proof.} As the proof is similar to those in \cite{18}, so we omit it here.

Define the following  BSDEs as
\begin{eqnarray}
dP_{t}
\hspace{-3mm}&=&\hspace{-3mm}-[Q_{t}+P_{t}A_{t}+A'_{t}P_{t}+C'_{t}P_{t}C_{t}
+\bar{P}_{t}C_{t}+C'_{t}\bar{P}_{t}]dt+\bar{P}_{t}dw(t), \label{f4}\\
dM(t,\theta-t)
\hspace{-3mm}&=&\hspace{-3mm}-[M(t,\theta-t)A_{t}+\bar{M}(t,\theta-t)C_{t}]dt
+\bar{M}(t,\theta-t)dw(t), \nonumber\\
\hspace{-3mm}&&\hspace{-3mm}\ \ t\in[0,T],\ \ \theta\in[t,T], \label{f5}\\
M(t,0)\hspace{-3mm}&=&\hspace{-3mm} B_{t}'P_{t}+D_{t}'\bar{P}_{t}+D_{t}'P_{t}C_{t}
-\int^{T}_{t}\Big\{E[M(t,\theta-t)B_{t}+\bar{M}(t,\theta-t)D_{t}]'\nonumber\\
\hspace{-3mm}&&\hspace{-3mm}\times \Upsilon^{-1}_{\theta}M(t,\theta-t)\Big\}d\theta, \label{f8}
\end{eqnarray}
with terminal values $P_{T}$. The above BSDEs (\ref{f4})-(\ref{f8}) will be called the coupled stochastic Riccati-type differential equations for discussing.

In view of the above preliminaries, the solution  to Problem 1 will be given in the sequence.
\begin{theorem}
  Let Assumption 1 be satisfied and the BSDEs (\ref{f4})-(\ref{f8}) has a solution. Problem 1 is uniquely solvable if and only if $\Upsilon_{t}>0$ for $t\in[0,T]$, where
\begin{eqnarray}
\Upsilon_{t}=E[R_{t}+D'_{t}P_{t}D_{t}]. \label{f6}
\end{eqnarray}
In this case, the optimal controller is
\begin{eqnarray}
u(t)=-\Upsilon^{-1}_{t}E[M(t,0)x(t)], \label{f7}
\end{eqnarray}
and the optimal cost function can be expressed as
\begin{eqnarray}
J^{\ast}_{T}=E\{x'(0)P_{0}x(0)
-\int^{T}_{0}(M(0,\theta)x(0))'
\Upsilon^{-1}_{\theta}E[M(0,\theta)x(0)]d\theta\},\label{f9}
\end{eqnarray}
moreover, the relationship between costate $p(t)$ in (\ref{f3}) and state $x(t)$ can be presented as
\begin{eqnarray}
p(t)&=&P_{t}x(t)
-\int^{T}_{t}\Big\{(M(t,\theta-t))'
\Upsilon^{-1}_{\theta}E[M(t,\theta-t)x(t)]\Big\}d\theta.\label{f10}
\end{eqnarray}

\end{theorem}
\emph{Proof.} \emph{Sufficiency:} \ We will verify that when $\Upsilon_{t}>0$, for any $t\in[0,T]$,  Problem 1 has a unique solution. \\
For convenience,  let
\begin{eqnarray}
dP_{t}&=&F_{t}dt+H_{t}dw(t),\label{f11}\\
d M(t,\theta-t)&=&\Sigma(t,\theta-t)dt+\Lambda(t,\theta-t)dw(t),\label{f12}
\end{eqnarray}
in which
\begin{eqnarray}
F_{t}
&=&-[Q_{t}+P_{t}A_{t}+A'_{t}P_{t}+C'_{t}P_{t}C_{t}
+\bar{P}_{t}C_{t}+C'_{t}\bar{P}_{t}],\ \  H_{t}=\bar{P}_{t}, \\ \label{f13}
\Sigma(t,\theta-t)
&=&-[M(t,\theta-t)A_{t}+\bar{M}(t,\theta-t)C_{t}],
\ \ \Lambda(t,\theta-t)=\bar{M}(t,\theta-t). \label{f14}
\end{eqnarray}\\
Hence, $d\Big[P_{t}x(t)
-\int^{T}_{t}\Big\{(M(t,\theta-t))'
\Upsilon^{-1}_{\theta}E[M(t,\theta-t)x(t)]\Big\}d\theta \Big]$ can be calculated by
\begin{eqnarray}
&&d\Big[P_{t}x(t)
-\int^{T}_{t}\Big\{(M(t,\theta-t))'
\Upsilon^{-1}_{\theta}E[M(t,\theta-t)x(t)]\Big\}d\theta \Big]\nonumber\\
&=&dP_{t}\cdot x(t)+P_{t}\cdot dx_{t}+dP_{t}\cdot dx_{t}
-\int^{T}_{t}\Big\{(d M(t,\theta-t))'
\Upsilon^{-1}_{\theta}E[M(t,\theta-t)x(t)]\nonumber\\
&&+(M(t,\theta-t))'
\Upsilon^{-1}_{\theta}E[d M(t,\theta-t)x(t)+M(t,\theta-t)dx_{t}
+ d M(t,\theta-t)dx_{t}]\Big\}d\theta\nonumber\\
&&+M'(t,0)\Upsilon^{-1}_{t}E[M(t,0)x(t)]dt\nonumber\\
&=&[F_{t}dt+H_{t}dw(t)]\cdot x(t)+P_{t}\cdot [(A_{t}x(t)+B_{t}u(t))dt+(C_{t}x(t)+D_{t}u(t))dw(t)]\nonumber\\
&&+H_{t}(C_{t}x(t)+D_{t}u(t))dt
-\int^{T}_{t}\Big\{(\Sigma(t,\theta-t)dt+\Lambda(t,\theta-t)dw(t))'
\Upsilon^{-1}_{\theta}E[M(t,\theta-t)x(t)]\nonumber\\
&&+(M(t,\theta-t))'
\Upsilon^{-1}_{\theta}E\Big[(\Sigma(t,\theta-t)dt+\Lambda(t,\theta-t)dw(t))x(t)
+M(t,\theta-t)[(A_{t}x(t)+B_{t}u(t))dt\nonumber\\
&&+(C_{t}x(t)+D_{t}u(t))dw(t)]
+\Lambda(t,\theta-t)(C_{t}x(t)+D_{t}u(t))dt\Big]\Big\}d\theta
+M'(t,0)\Upsilon^{-1}_{t}E[M(t,0)x(t)]dt\nonumber\\
&=&\Big\{\Big(F_{t}+P_{t}A_{t}+ H_{t}C_{t}\Big)x(t)+\Big(P_{t}B_{t}+H_{t}D_{t}-\int^{T}_{t}(M(t,\theta-t))'
\Upsilon^{-1}_{\theta}E[M(t,\theta-t)B_{t}
\nonumber\\
&&+\Lambda(t,\theta-t)D_{t}]d\theta\Big)u(t)-\int^{T}_{t}\Big((\Sigma(t,\theta-t))'
\Upsilon^{-1}_{\theta}E[M(t,\theta-t)x(t)]+(M(t,\theta-t))'
\Upsilon^{-1}_{\theta}\nonumber\\
&&\times E[\Sigma(t,\theta-t)x(t)
+M(t,\theta-t)A_{t}x(t)+\Lambda(t,\theta-t)C_{t}x(t)]\Big)d\theta
+M'(t,0)\Upsilon^{-1}_{t}E[M(t,0)x(t)]\Big\}dt \nonumber\\ &&+\Big\{H_{t}x(t)+P_{t}C_{t}x(t)+P_{t}D_{t}u(t)
-\int^{T}_{t}(\Lambda(t,\theta-t))'
\Upsilon^{-1}_{\theta}E[M(t,\theta-t)x(t)]d\theta\Big\}dw(t).\label{f15}
\end{eqnarray}
On this basis, applying It$\hat{o}$'s formula to $x'(t)\Big[P_{t}x(t)
-\int^{T}_{t}\Big\{(M(t,\theta-t))'
\Upsilon^{-1}_{\theta}E[M(t,\theta-t)x(t)]\Big\}d\theta \Big]$, taking integral from 0 to $T$ and then expectation,  we have that
\begin{eqnarray}
&&E\int^{T}_{0}d\Big\{x'(t)\Big[P_{t}x(t)
-\int^{T}_{t}(M(t,\theta-t))'
\Upsilon^{-1}_{\theta}E[M(t,\theta-t)x(t)]d\theta \Big]\Big\}\nonumber\\
&=&E\int^{T}_{0}\Big\{dx'(t)\cdot\Big[P_{t}x(t)
-\int^{T}_{t}(M(t,\theta-t))'
\Upsilon^{-1}_{\theta}E[M(t,\theta-t)x(t)]d\theta \Big]\nonumber\\
&&+x'(t)\cdot d\Big[P_{t}x(t)
-\int^{T}_{t}(M(t,\theta-t))'
\Upsilon^{-1}_{\theta}E[M(t,\theta-t)x(t)]d\theta \Big]\nonumber\\
&&+dx'(t)\cdot d\Big[P_{t}x(t)
-\int^{T}_{t}(M(t,\theta-t))'
\Upsilon^{-1}_{\theta}E[M(t,\theta-t)x(t)]d\theta \Big]\Big\}\nonumber
\end{eqnarray}
\begin{eqnarray}
&=&E\int^{T}_{0}\Big\{x'(t)[A'_{t}P_{t}+F_{t}+P_{t}A_{t}+H_{t}C_{t}
+C'_{t}H_{t}+C'_{t}P_{t}C_{t}]x(t)+u'(t)D'_{t}P_{t}D_{t}u(t)\nonumber\\
&&+x'(t)\Big(P_{t}B_{t}+H_{t}D_{t}+C'_{t}P_{t}D_{t}-\int^{T}_{t}(M(t,\theta-t))'
\Upsilon^{-1}_{\theta}E[M(t,\theta-t)B_{t}
+\Lambda(t,\theta-t)D_{t}]d\theta\Big)'u(t)\nonumber\\
&&-x'(t)\int^{T}_{t}(A_{t}M(t,\theta-t)+\Sigma(t,\theta-t)+C_{t}\Lambda(t,\theta-t))
\Upsilon^{-1}_{\theta}E[M(t,\theta-t)x(t)]d\theta\nonumber\\
&&
+u'(t)\Big[(B'_{t}P_{t}+D'_{t}P_{t}C_{t}+D'_{t}H_{t})x(t)-D'_{t}\int^{T}_{t}(\Lambda(t,\theta-t))'
\Upsilon^{-1}_{\theta}E[M(t,\theta-t)x(t)]d\theta\nonumber\\
&&-B'_{t}\int^{T}_{t}(M(t,\theta-t))'
\Upsilon^{-1}_{\theta}E[M(t,\theta-t)x(t)]d\theta\Big]
+x'(t)(M(t,0))'\Upsilon^{-1}_{t}E[M(t,0)x(t)]\nonumber\\
&&-x'(t)\int^{T}_{t}(M(t,\theta-t))'
\Upsilon^{-1}_{\theta} E[(\Sigma(t,\theta-t)+M(t,\theta-t)A_{t}+\Lambda(t,\theta-t)C_{t})x(t)]d\theta
\Big\}.\label{f16}
\end{eqnarray}
In view of (\ref{f11})- (\ref{f14}) and (\ref{f16}), it yields that
\begin{eqnarray}
J_{T}&=&E\Big\{x'(0)P_{0}x(0)
-\int^{T}_{0}(M(0,\theta)x(0))'
\Upsilon^{-1}_{\theta}E[M(0,\theta)x(0)]d\theta\nonumber\\
&&+\int^{T}_{0}\Big(u(t)
+\Upsilon^{-1}_{t}E[M(t,0)x(t)]\Big)'\Upsilon_{t}\Big(u(t)
+\Upsilon^{-1}_{t}E[M(t,0)x(t)]\Big)dt\Big\}.\label{f17}
\end{eqnarray}
Since $\Upsilon_{t}>0$, the optimal controller can be obtained, i.e.,
\begin{eqnarray}
u(t)
=-\Upsilon^{-1}_{t}E[M(t,0)x(t)],\label{f18}
\end{eqnarray}
 and the optimal cost function can also be presented as
 \begin{eqnarray}
J^{\ast}_{T}=E\Big\{x'(0)P_{0}x(0)
-\int^{T}_{0}(M(0,\theta)x(0))'
\Upsilon^{-1}_{\theta}E[M(0,\theta)x(0)]d\theta\Big\}.\label{f19}
\end{eqnarray}
\emph{Necessity:} Next we will illustrate the result that if Problem 1 has a unique solution, $\Upsilon_{t}>0$ will be satisfied.\\
%The outline of the proof can be expressed as: Assume that (\ref{f10}) is a solution to FBSDEs (\ref{f3}) by applying It$\hat{o}$ formula to (\ref{f10}) and then compared with the expression in (\ref{f3}). Finally, in view of Lemma 1, the solution to FBSDEs (\ref{f3}) can be derived.\\
From Lemma 1, when Problem 1 is solvable, the FBSDEs (\ref{f3}) exists a solution.
Assume that the solution to  FBSDEs (\ref{f3}) can be expressed as
\begin{eqnarray}
p(t)&=&P_{t}x(t)
-\int^{T}_{t}\Big\{(M(t,\theta-t))'
\Upsilon^{-1}_{\theta}E[M(t,\theta-t)x(t)]\Big\}d\theta,\label{f21}\\
q(t)&=&\bar{P}_{t}x(t)+P_{t}C_{t}x(t)+P_{t}D_{t}u(t)
-\int^{T}_{t}(\bar{M}(t,\theta-t))'
\Upsilon^{-1}_{\theta}E[M(t,\theta-t)x(t)]d\theta,\label{f22}
\end{eqnarray}
where $(P_{t},\bar{P}_{t})$ and $(M(t,\theta-t),\bar{M}(t,\theta-t))$ satisfy the BSDEs (\ref{f4})-(\ref{f8}).\\
In the sequence, we will verify that $dp(t)+[A_{t}'p(t)+C_{t}'q(t)+Q_{t}x(t)]dt-q(t)dw(t)=0$. Applying It$\hat{o}$'s formula to $p(t)$ in (\ref{f21}), similar to the lines of (\ref{f15}), we can obtain that
\begin{eqnarray}
&&dp(t)+[A_{t}'p(t)+C_{t}'q(t)+Q_{t}x(t)]dt-q(t)dw(t)\nonumber\\
&=&\Big\{\Big(-[Q_{t}+P_{t}A_{t}+A'_{t}P_{t}+C'_{t}P_{t}C_{t}
+\bar{P}_{t}C_{t}+C'_{t}\bar{P}_{t}]+P_{t}A_{t}+ \bar{P}_{t}C_{t}\Big)x(t)\nonumber\\
&&+\Big(P_{t}B_{t}+\bar{P}_{t}D_{t}-\int^{T}_{t}(M(t,\theta-t))'
\Upsilon^{-1}_{\theta}E[M(t,\theta-t)B_{t}
+\bar{M}(t,\theta-t)D_{t}]d\theta\Big)u(t)\nonumber\\
&&-\int^{T}_{t}\Big(-[M(t,\theta-t)A_{t}+\bar{M}(t,\theta-t)C_{t}])'
\Upsilon^{-1}_{\theta}E[M(t,\theta-t)x(t)]\nonumber\\
&&+(M(t,\theta-t))'
\Upsilon^{-1}_{\theta}E[-(M(t,\theta-t)A_{t}+\bar{M}(t,\theta-t)C_{t})x(t)
+M(t,\theta-t)A_{t}x(t)\nonumber\\
&&+\bar{M}(t,\theta-t)C_{t}x(t)]\Big)d\theta
+M'(t,0)\Upsilon^{-1}_{t}E[M(t,0)x(t)]\Big\}dt \nonumber\\ &&+\Big\{\bar{P}_{t}x(t)+P_{t}C_{t}x(t)+P_{t}D_{t}u(t)
-\int^{T}_{t}(\bar{M}(t,\theta-t))'
\Upsilon^{-1}_{\theta}E[M(t,\theta-t)x(t)]d\theta\Big\}dw(t)\nonumber\\
&&-q(t)dw(t)+\Big\{A_{t}'\Big[P_{t}x(t)
-\int^{T}_{t}\Big\{(M(t,\theta-t))'
\Upsilon^{-1}_{\theta}E[M(t,\theta-t)x(t)]\Big\}d\theta\Big]+Q_{t}x(t)\nonumber\\
&&+C_{t}'\Big[\bar{P}_{t}x(t)+P_{t}C_{t}x(t)+P_{t}D_{t}u(t)
-\int^{T}_{t}(\bar{M}(t,\theta-t))'
\Upsilon^{-1}_{\theta}E[M(t,\theta-t)x(t)]d\theta\Big]
\Big\}dt. \label{f24}
\end{eqnarray}
Plugging (\ref{f21}) and (\ref{f22}) into (\ref{f3}), we have that
\begin{eqnarray}
0&=&E\Big\{B_{t}'[P_{t}x_{t}-\int^{T}_{t}\Big\{(M(t,\theta-t))'
\Upsilon^{-1}_{\theta}E[M(t,\theta-t)x_{t}]\Big\}d\theta]+R_{t}u(t)\nonumber\\
&&+D_{t}'[\bar{P}_{t}x_{t}+P_{t}C_{t}x(t)+P_{t}D_{t}u(t)
-\int^{T}_{t}(\bar{M}(t,\theta-t))'
\Upsilon^{-1}_{\theta}E[M(t,\theta-t)x_{t}]d\theta]\Big\}\nonumber\\
&=&E[R_{t}+D_{t}'P_{t}D_{t}]u(t)+E\Big\{\Big[B_{t}'P_{t}+D_{t}'\bar{P}_{t}
+D_{t}'P_{t}C_{t}\nonumber\\
&&-\int^{T}_{t}E[M(t,\theta-t)B_{t}+\bar{M}(t,\theta-t)D_{t}]'
\Upsilon^{-1}_{\theta}M(t,\theta-t)d\theta\Big]x(t)\Big\}\nonumber\\
&=&\Upsilon_{t}u(t)+E[M(t,0)x(t)],\label{f25}
\end{eqnarray}
 it yields that
\begin{eqnarray}
u(t)=-\Upsilon^{-1}_{t}E[M(t,0)x(t)].\label{f26}
\end{eqnarray}
In view of (\ref{f24}) and (\ref{f26}), we can obtain that $dp(t)+[A_{t}'p(t)+C_{t}'q(t)+Q_{t}x(t)]dt-q(t)dw(t)=0$. Noting Lemma 1, i.e., when Problem 1 is uniquely solvable, the FBSDEs (\ref{f3}) has a unique solution, we know that (\ref{f21}) and (\ref{f22}) is the unique solution of FBSDEs (\ref{f3}).\\
Applying It$\hat{o}$'s formula to $x'(t)p(t)$, then similar to the line of (\ref{f16}), we obtain that
\begin{eqnarray}
J_{T}&=&E\Big\{x'(0)P_{0}x(0)
-\int^{T}_{0}(M(0,\theta)x(0))'
\Upsilon^{-1}_{\theta}E[M(0,\theta)x(0)]d\theta\nonumber\\
&&+\int^{T}_{0}\Big(u(t)
+\Upsilon^{-1}_{t}E[M(t,0)x(t)]\Big)'\Upsilon_{t}\Big(u(t)
+\Upsilon^{-1}_{t}E[M(t,0)x(t)]\Big)dt\Big\}.\label{f27}
\end{eqnarray}
In what follows, $\Upsilon_{t}>0$ will be proved.\\
Let $\lambda_{t}$ be any fixed eigenvalue of the matrix $\Upsilon_{t}, t\in[0, T]$. We will show that $mes(\{t\in[0, T]|\lambda_{t}<0\})=0$, where $mes$ denotes the Lebesgue measure. Let $\nu_{\lambda}(t)$ be a unit eigenvector (i.e., $\nu_{\lambda}^{'}(t)\nu_{\lambda}(t)=1$) associated with the eigenvalue $\lambda_{t}$. Define $I_{n}$ as the indicator function of the set $\{t\in[0, T]|\lambda_{t}<-\frac{1}{n}\}, n=1,2,\cdots.$ Fix a scalar $\delta\in \mathcal{R}$ and consider the state trajectory $x(\cdot)$ of system (\ref{f1}) under the feedback control
\begin{eqnarray}
u(t)=
\left\{
\begin{array}{lll}
0, \ \ \ \ \ \ \ \ \ \ \ \ \ \ \ \ \ \ \ \ if \ \lambda_{t}=0 \\
\frac{\delta I_{n}(t)}{|\lambda_{t}|^{\frac{1}{2}}}\nu_{\lambda_{t}}
-\Upsilon_{t}^{\dagger}E[M(t,0)x(t)], \ if \ \lambda_{t}\neq0
\end{array}
\right.     \label{f28}
\end{eqnarray}
On this basis, using the controller (\ref{f28}), then $J_{T}$ in (\ref{f27}) can be further calculated as
\begin{eqnarray}
J_{T}&=&E\Big\{x'(0)P_{0}x(0)
-\int^{T}_{0}(M(0,\theta))'
\Upsilon^{-1}_{\theta}E[M(0,\theta)x(0)]d\theta
-\delta^{2}\int_{I_{n}}dt\Big\}\nonumber\\
&=&\frac{1}{2}E\Big\{x'(0)P_{0}x(0)
-\int^{T}_{0}(M(0,\theta))'
\Upsilon^{-1}_{\theta}E[M(0,\theta)x(0)]d\theta\nonumber\\
&&-\delta^{2}mes\Big(\{t\in[0,T]|\lambda_{t}<-\frac{1}{n}\}\Big)\Big\}.\label{f29}
\end{eqnarray}
 If $mes\Big(\{t\in[0,T]|\lambda_{t}<-\frac{1}{n}\}\Big)>0$, then by letting $\delta\rightarrow\infty$, we obtain $J_{T}\rightarrow-\infty$, which contradicts with $J_{T}\geq 0$. Hence $mes\Big(\{t\in[0,T]|\lambda_{t}<-\frac{1}{n}\}\Big)=0$.
 Since $\{t\in[0,T]|\lambda_{t}<0\}=\bigcup^{\infty}_{n=1}\{t\in[0,T]|\lambda_{t}<-\frac{1}{n}\}$, we conclude that $mes\Big(\{t\in[0,T]|\lambda_{t}<0\}\Big)=0$, i.e., $\Upsilon_{t}\geq 0$. Finally, by virtue of the unique solvability of Problem 1 and (\ref{f27}), the positive-definiteness of $\Upsilon_{t}$ follows. This completes the proof.
 \begin{assumption}
 $A_{t},B_{t},C_{t},D_{t}, Q_{t}, R_{t}$ are deterministic matrix-valued functions with suitable sizes, and $Q_{t}, R_{t}, P_{T}$ are positive semi-definite.
 \end{assumption}
\begin{corollary}
Let Assumption 2 be satisfied. Problem 1 is uniquely solvable if and only if
\begin{eqnarray}
\tilde{\Upsilon}_{t}=R_{t}+D_{t}'P_{1}(t)D_{t}, \label{f30}
\end{eqnarray}
is strictly positive, in which
\begin{eqnarray}
\dot{P}_{1}(t)+Q_{t}+P_{1}(t)A_{t}+A_{t}'P_{1}(t)+C_{t}'P_{1}(t)C_{t}=0, P_{1}(T)=P_{T}. \label{f31}
\end{eqnarray}
In this case,  the  optimal control is
\begin{eqnarray}
u(t)=-\tilde{\Upsilon}^{-1}_{t}[B_{t}P_{2}(t)+D_{t}'P_{1}(t)C_{t}]Ex_{t}, \label{f34}
\end{eqnarray}
where
\begin{eqnarray}
\dot{P}_{2}(t)+P_{2}(t)\tilde{A}_{t}+\tilde{A}_{t}'P_{2}(t)+\tilde{Q}_{t}
-P_{2}(t)B_{t}'\tilde{\Upsilon}^{-1}_{t}B_{t}P_{2}(t)=0,P_{2}(T)=P_{T}\label{f039}
\end{eqnarray}
with $\tilde{A}_{t}:=A_{t}-B_{t}'\tilde{\Upsilon}^{-1}_{t}D_{t}'P_{1}(t)C_{t}$,  $\tilde{Q}_{t}:=Q_{t}+C_{t}'P_{1}(t)C_{t}
-C_{t}'P_{1}(t)D_{t}\tilde{\Upsilon}^{-1}_{t}D_{t}'P_{1}(t)C_{t}$,
and the optimal cost can be obtained that
\begin{eqnarray}
J_{T}= E\Big\{x(0)'P_{1}(0)[x(0)-Ex(0)]+x(0)'P_{2}(0) Ex(0)\Big\}.\label{f35}
\end{eqnarray}
Moreover, the relationship between the costate $p(t)$ and state $x(t)$ is
\begin{eqnarray}
p(t)=P_{1}(t)[x(t)-Ex(t)]+P_{2}(t)Ex(t).\label{f36}
\end{eqnarray}
\end{corollary}
\emph{proof.} \emph{``Sufficiency:"} When $\tilde{\Upsilon}_{t}$ in (\ref{f30}) is positive, we will illustrate the unique solvability of Problem 1 under Assumption 2.
Applying It$\hat{o}$ formula to $P_{1}(t)[x(t)-Ex(t)]+P_{2}(t)Ex(t)$, we have that
\begin{eqnarray}
&&d\Big\{P_{1}(t)[x(t)-Ex(t)]+P_{2}(t)Ex(t)\Big\}\nonumber\\
&=&\dot{P}_{1}(t)[x(t)-Ex(t)]dt+P_{1}(t)d[x(t)-Ex(t)]+\dot{P}_{2}(t)Ex(t)dt
+P_{2}(t)dEx(t)\nonumber\\
&=&\dot{P}_{1}(t)[x(t)-Ex(t)]dt+P_{1}(t)[(A_{t}x(t)+B_{t}u(t))dt
+(C_{t}x(t)+D_{t}u(t))dw(t)\nonumber\\
&&-(A_{t}Ex(t)+B_{t}u(t))dt]+\dot{P}_{2}(t)Ex(t)dt
+P_{2}(t)(A_{t}Ex(t)+B_{t}u(t))dt\nonumber\\
&=&\Big\{[\dot{P}_{1}(t)+P_{1}(t)A_{t}]x(t)+[-\dot{P}_{1}(t)
-P_{1}(t)A_{t}+\dot{P}_{2}(t)+P_{2}(t)A_{t}]Ex(t)+P_{2}(t)B_{t}u(t)\Big\}dt\nonumber\\
&&+[P_{1}(t)C_{t}x(t)+P_{1}(t)D_{t}u(t)]dw(t).\label{f37}
\end{eqnarray}
In view of (\ref{f37}), it yields that
\begin{eqnarray}
&&d\Big\{x(t)'\left[P_{1}(t)(x(t)-Ex(t))+P_{2}(t)Ex(t)\right]\Big\}\nonumber\\
&=&\Big\{(A_{t}x(t)+B_{t}u(t))'\left[P_{1}(t)(x(t)-Ex(t))+P_{2}(t)Ex(t)\right]
+x(t)'[\dot{P}_{1}(t)+P_{1}(t)A_{t}]x(t)\nonumber\\
&&+x(t)'[-\dot{P}_{1}(t)
-P_{1}(t)A_{t}+\dot{P}_{2}(t)+P_{2}(t)A_{t}]Ex(t)+x(t)'P_{2}(t)B_{t}u(t)\nonumber\\
&&+(C_{t}x(t)+D_{t}u(t))'[P_{1}(t)C_{t}x(t)+P_{1}(t)D_{t}u(t)]\Big\}dt
+\Big\{(C_{t}x(t)+D_{t}u(t))'\nonumber\\
&&\times\Big[P_{1}(t)[x(t)-Ex(t)]+P_{2}(t)Ex(t)\Big]
+x(t)'[P_{1}(t)C_{t}x(t)+P_{1}(t)D_{t}u(t)]\Big\}dw(t)\nonumber\\
&=&\Big\{x(t)'[\dot{P}_{1}(t)+P_{1}(t)A_{t}+A_{t}'P_{1}(t)+C_{t}'P_{1}(t)C_{t}]x(t)
+u(t)'D_{t}'P_{1}(t)D_{t}u(t)\nonumber\\
&&+x(t)'[-A_{t}'P_{1}(t)+A_{t}'P_{2}(t)
-\dot{P}_{1}(t)-P_{1}(t)A_{t}+\dot{P}_{2}(t)+P_{2}(t)A_{t}]Ex(t)\nonumber\\
&&+u(t)'[B_{t}'P_{1}(t)+D_{t}'P_{1}(t)C_{t}]x(t)
+x(t)'[P_{2}(t)B_{t}+C_{t}'P_{1}(t)D_{t}]u(t)\nonumber\\
&&+u(t)'[-B_{t}'P_{1}(t)+B_{t}'P_{2}(t)]Ex(t)\Big\}dt+\Big\{(C_{t}x(t)+D_{t}u(t))'
[P_{1}(t)(x(t)-Ex(t))\nonumber\\
&&+P_{2}(t)Ex(t)]+x(t)'[P_{1}(t)C_{t}x(t)+P_{1}(t)D_{t}u(t)]\Big\}dw(t).\label{f38}
\end{eqnarray}
Taking the integral from 0 to $T$ and the expectation on both sides of (\ref{f38}), it yields that
\begin{eqnarray}
J_{T}&=&E\Big[x(0)'\left[P_{1}(0)(x(0)-Ex(0))+P_{2}(0)Ex(0)\right]\Big]\nonumber\\
&&+\int^{T}_{0}\Big\{Ex(t)'[\dot{P}_{2}(t)+P_{2}(t)A_{t}+A_{t}'P_{2}(t)
+C_{t}'P_{1}(t)C_{t}+Q_{t}]Ex(t)\nonumber\\
&&+u(t)'[B_{t}'P_{1}(t)+D_{t}'P_{1}(t)C_{t}]Ex(t)
+Ex(t)'[P_{2}(t)B_{t}+C_{t}'P_{1}(t)D_{t}]u(t)\nonumber\\
&&+u(t)'[-B_{t}'P_{1}(t)+B_{t}'P_{2}(t)]Ex(t)\Big\}dt\nonumber\\
&=&E\Big[x(0)'\left[P_{1}(0)(x(0)-Ex(0))+P_{2}(0)Ex(0)\right]\Big]\nonumber\\
&&+\int^{T}_{0}\Big\{Ex(t)'(B_{t}P_{2}(t)+D_{t}'P_{1}(t)C_{t})'
\tilde{\Upsilon}^{-1}_{t}(B_{t}P_{2}(t)+D_{t}'P_{1}(t)C_{t})Ex(t)\nonumber\\
&&+u(t)'[B_{t}'P_{2}(t)+D_{t}'P_{1}(t)C_{t}]Ex(t)
+Ex(t)'[P_{2}(t)B_{t}+C_{t}'P_{1}(t)D_{t}]u(t)\Big\}dt\nonumber\\
&=&E\Big[x(0)'\left[P_{1}(0)(x(0)-Ex(0))+P_{2}(0)Ex(0)\right]\Big]
+\int^{T}_{0}[u(t)+\tilde{\Upsilon}^{-1}_{t}(B_{t}P_{2}(t)+D_{t}'P_{1}(t)C_{t})Ex(t)]'\nonumber\\
&&\times
\tilde{\Upsilon}_{t}[u(t)+\tilde{\Upsilon}^{-1}_{t}(B_{t}P_{2}(t)+D_{t}'P_{1}(t)C_{t})Ex(t)]dt.\label{f39}
\end{eqnarray}
From (\ref{f39}) and $\tilde{\Upsilon}_{t}>0$, we can obtain that the optimal control is
\begin{eqnarray}
u(t)=-\tilde{\Upsilon}^{-1}_{t}(B_{t}P_{2}(t)+D_{t}'P_{1}(t)C_{t})Ex(t),\label{f40}
\end{eqnarray}
and the optimal cost is
\begin{eqnarray}
J_{T}=E\Big[x(0)'\left[P_{1}(0)(x(0)-Ex(0))+P_{2}(0)Ex(0)\right]\Big].\label{f41}
\end{eqnarray}
\emph{``Necessary:"} Under Assumption 2, if Problem 1 has a unique solution, $\tilde{\Upsilon}_{t}>0$ will be proved. From (\ref{f37}), similar to the lines of the necessary part in Theorem 1, we can derive the proof, so we omit it here.

\begin{remark}
 The above result presented in Corollary 1 is parallel to the result for the case of $B^{2}_{t} = 0$ and $D^{2}_{t} = 0$  in \cite{15}. Concretely, under Assumption 2, (\ref{f4})-(\ref{f8}) can be rewritten as
\begin{eqnarray}
0&=&\dot{P}_{t}+Q_{t}+P_{t}A_{t}+A_{t}'P_{t}+C_{t}'P_{t}C_{t}, \label{f42}\\
\frac{\partial M(t,\theta-t)}{\partial t}&=&-M(t,\theta-t)A_{t}, \label{f43}\\
M(t,0)&=&B_{t}'\Big[P_{t}-\int^{T}_{t}M(t,\theta-t)'
\Upsilon^{-1}_{\theta}M(t,\theta-t)d\theta\Big]+D_{t}'P_{t}C_{t}, \label{f46}\\
\Upsilon_{t}&=&R_{t}+D_{t}'P_{t}D_{t}, \label{f44}\\
u(t)&=&-\Upsilon^{-1}_{t}M(t,0)Ex_{t}. \label{f45}
\end{eqnarray}
Further, let
\begin{eqnarray}
F_{t}=-\int^{T}_{t}M(t,\theta-t)'\Upsilon^{-1}_{\theta}M(t,\theta-t)d\theta,\label{f47}
\end{eqnarray}
and its derivative can be obtained that
\begin{eqnarray}
\dot{F}_{t}=-F_{t}A_{t}-A_{t}'F_{t}+M(t,0)'\Upsilon^{-1}_{t}M(t,0).\label{f48}
\end{eqnarray}
Let $P_{2}(t)=P_{t}+F_{t}$, from (\ref{f42}) and (\ref{f48}), it yields that
\begin{eqnarray}
\dot{P}_{2}(t)&=&\dot{P}_{t}+\dot{F}_{t}\nonumber\\
&=&-(P_{t}+F_{t})A_{t}-A_{t}'(P_{t}+F_{t})-C_{t}'P_{t}C_{t}-Q_{t}
+M(t,0)'\Upsilon^{-1}_{t}M(t,0)\nonumber\\
&=&-P_{2}(t)A_{t}-A_{t}'P_{2}(t)-C_{t}'P_{t}C_{t}-Q_{t}
+(B_{t}'P_{2}(t)+D_{t}'P_{t}C_{t})'\Upsilon^{-1}_{t}(B_{t}'P_{2}(t)\nonumber\\
&&+D_{t}'P_{t}C_{t}).\label{f49}
\end{eqnarray}
Therefore, (\ref{f10}) can be reexpressed as
\begin{eqnarray}
p(t)=P_{t}[x(t)-Ex(t)]+P_{2}(t)Ex(t).\label{f50}
\end{eqnarray}
In view of these,  we can obtain that the main result presented in Theorem 1 can be reduced to the result in Corollary 1, i.e., the case of $B^{2}_{t} = 0$ and $D^{2}_{t} = 0$  in \cite{15}.
\end{remark}

\section{Conclusion}
In this paper, we have solved the deterministic LQ control of It$\hat{o}$ stochastic systems with random coefficients by presenting the necessary and sufficient solving conditions and explicit controller via a new coupled Riccati-type equation.

As compared with the LQ control problem with deterministic coefficients or $\mathcal{F}_{t}$-adapted control, the considered problem in this paper is very involved.  The key technique  for us to derive the results are the analytical solution to the FBSDEs  originated from maximum principle.

It can be shown from this paper and our earlier works \cite{16} and \cite{17}  that the technique of solving FBSDEs  is very powerful to optimal control.  Actually, it has been successfully applied  to solve many difficult problems such as stochastic control with delay, irregular LQ control, LQ control with asymmetric information and LQ control in networked control systems and so on; see  \cite{19}-\cite{23} and references therein.

\section*{Acknowledgements}

 The authors would like to thank Prof. Shanjian Tang for his valuable discussions.

\bibliographystyle{plain}        % Include this if you use bibtex
%\bibliography{autosam}           % and a bib file to produce the
                                 % bibliography (preferred). The
                                 % correct style is generated by
                                 % Elsevier at the time of printing.

\end{document}